\documentclass{agtart_a}
\pdfoutput=1

\usepackage{graphicx}
\usepackage{color}
\usepackage[all]{xy}

\title{Totally geodesic surfaces and homology}

\author{Jason DeBlois}
\givenname{Jason}
\surname{DeBlois}
\address{Department of Mathematics\\
 The University of Texas at Austin\\\newline
  1 University Station C1200\\ 
  Austin, TX 78712-0257\\USA}
\email{jdeblois@math.utexas.edu}
\urladdr{}

\volumenumber{6}
\issuenumber{}
\publicationyear{2006}
\papernumber{51}
\startpage{1413}
\endpage{1428}

\doi{}
\MR{}
\Zbl{}


\received{15 February 2006}
\revised{5 July 2006}
\accepted{7 July 2006}
\published{25 September 2006}
\publishedonline{25 September 2006}
\proposed{}
\seconded{}
\corresponding{}
\editor{}
\version{}

\arxivreference{math.GT/0601561}

\makeatletter
\def\cnewtheorem#1[#2]#3{\newtheorem{#1}{#3}[section]
\expandafter\let\csname c@#1\endcsname\c@thm}
\makeatother

\let\xysavmatrix\xymatrix
\def\xymatrix{\disablesubscriptcorrection\xysavmatrix}
\AtBeginDocument{\let\bar\wbar\let\tilde\wtilde\let\hat\what}
\def\SetFigFont#1#2#3#4#5{\small}%
\def\adjustlabel<#1,#2>#3{\smash{\rlap{\kern #1 \raise #2\hbox{#3}}}}
\def\tphi{{}\mskip2.5mu\tilde{\mskip-2.5mu \vphantom{t}\smash{\phi}\mskip-1mu}\mskip1mu}
\makeop{tr}

\theoremstyle{plain}
\newtheorem{thm}{Theorem}
\newtheorem{lem}{Lemma}
\newtheorem*{thm1}{Theorem}

\newtheorem*{fact}{Fact}
\newtheorem*{cor1}{Corollary}
\newtheorem*{conj}{Conjecture}
\newtheorem*{ques}{Question}

\theoremstyle{remark}

\newtheorem*{remark}{Remark}

\begin{document}

\begin{asciiabstract}
We construct examples of hyperbolic rational homology spheres and
hyperbolic knot complements in rational homology spheres containing
closed embedded totally geodesic surfaces.
\end{asciiabstract}

\begin{abstract} 
We construct examples of hyperbolic rational homology spheres and
hyperbolic knot complements in rational homology spheres containing
closed embedded totally geodesic surfaces.
\end{abstract}

\maketitle

\section{Introduction}

Let $M=\mathbb{H}^3/ \Gamma$, $\Gamma \subset \mathrm{PSL}_2(\mathbb{C})$ be an orientable hyperbolic $3$--manifold, and let $f\co F \rightarrow M$ be a proper immersion of a connected, orientable surface of genus at least $2$ such that $f_{*}\co\pi_1(F) \rightarrow \Gamma$ is injective.
$F$ (or more precisely $(f,F)$) is said to be {\em  totally geodesic\/} if $f_{*}(\pi_{1}(F)) \subset \Gamma$ is conjugate into $\mathrm{PSL}_2(\mathbb{R})$.  Thurston and Bonahon have described the geometry of surface groups in hyperbolic $3$--manifolds as falling into three classes: doubly degenerate groups, quasi-Fuchsian groups and groups with accidental parabolics.  The class of totally geodesic surface groups is a ``positive codimension'' subclass of the quasi-Fuchsian groups, so one may expect that hyperbolic $3$--manifolds containing totally geodesic surface groups are special.

Indeed, the presence of a totally geodesic surface in a hyperbolic $3$--manifold has important topological implications.  Long showed that immersed totally geodesic surfaces lift to embedded nonseparating surfaces in finite covers \cite{Long}, proving the virtual Haken and virtually positive $\beta_1$ conjectures for hyperbolic manifolds containing totally geodesic surfaces.  Given this, it is natural to wonder about the extent to which topology constrains the existence of totally geodesic surfaces in hyperbolic $3$--manifolds.  Menasco--Reid have made the following conjecture \cite{MeR}:
\begin{conj}[Menasco--Reid] No hyperbolic knot complement in $\mathrm{S}^3$ contains a closed embedded totally geodesic surface.  
\end{conj}
They proved this conjecture for alternating knots.  The Menasco--Reid conjecture has been shown true for many other classes of knots, including almost alternating knots \cite{Aetal}, Montesinos knots \cite{Oer}, toroidally alternating knots \cite{Ad}, $3$--bridge and double torus knots \cite{IO} and knots of braid index 3 \cite{LP} and 4 \cite{Ma}.  For a knot in one of the above families, any closed essential surface in its complement has a topological feature which obstructs it from being even quasi-Fuchsian.  In general, however, one cannot hope to find such obstructions.  Adams--Reid have given examples of closed embedded quasi-Fuchsian surfaces in knot complements which volume calculations prove to be not totally geodesic \cite{AR}.

On the other hand, C\,Leininger has given evidence for a counterexample by constructing a sequence of hyperbolic knot complements in $\mathrm{S}^3$ containing closed embedded surfaces whose principal curvatures approach 0 \cite{Le}.  In this paper, we take an alternate approach to giving evidence for a counterexample.
\begin{thm} \label{cusped}
 There exist infinitely many hyperbolic knot complements in rational homology spheres containing closed embedded totally geodesic surfaces.  
\end{thm}
This answers a question of Reid---recorded as Question 6.2 in \cite{Le}---giving counterexamples to the natural generalization of the Menasco--Reid conjecture to knot complements in rational homology spheres.  Thus the conjecture, if true, must reflect a deeper topological feature of knot complements in $\mathrm{S}^3$ than simply their rational homology.

Prior to proving \fullref{cusped}, in \fullref{sec2} we prove the following theorem.
\begin{thm} \label{closed}
There exist infinitely many hyperbolic rational homology spheres containing closed embedded totally geodesic surfaces.
\end{thm}
This seems of interest in its own right, and the proof introduces many of the techniques we use in the proof of \fullref{cusped}.  Briefly, we find a two-cusped hyperbolic manifold containing an embedded totally geodesic surface which remains totally geodesic under certain orbifold surgeries on its boundary slopes and use the Alexander polynomial to show that branched covers of these surgeries have no rational homology.  In \fullref{sec3} we prove \fullref{cusped}.  In the final section, we give some idea of further directions and questions suggested by our approach.

\subsubsection*{Acknowledgements}

The author thanks Cameron Gordon, Richard Kent, Chris Lein\-inger, Jessica Purcell and Alan Reid for helpful conversations.  The author also thanks the Centre Interfacultaire Bernoulli at EPF Lausanne for their hospitality during part of this work.

\section[Theorem 2]{\fullref{closed}}\label{sec2}

Given a compact hyperbolic manifold $M$ with totally geodesic boundary of genus $g$, gluing it to its mirror image $\wwbar{M}$ along the boundary yields a closed manifold $DM$---the ``double'' of $M$---in which the former $\partial M$ becomes an embedded totally geodesic surface.  One limitation of this construction is that this surface contributes half of its first homology to the first homology of $DM$, so that $\beta_1(DM) \geq g$.  This is well known, but we include an argument to motivate our approach. Consider the relevant portion of the rational homology Mayer--Vietoris sequence for $DM$:  
\[ \cdots \rightarrow \mathrm{H}_1(\partial M,\mathbb{Q}) \stackrel{(i_*,-j_*)}{\longrightarrow} \mathrm{H}_1(M,\mathbb{Q}) \oplus \mathrm{H}_1(\wwbar{M},\mathbb{Q}) \rightarrow \mathrm{H}_1(DM,\mathbb{Q}) \rightarrow 0 \]
The labeled maps $i_*$ and $j_*$ are the maps induced by inclusion of the surface into $M$ and $\wwbar{M}$, respectively.  Recall that by the ``half lives, half dies'' lemma (see eg Hatcher \cite[Lemma 3.5]{Ha}), the dimension of the kernel of $i_*$ is equal to $g$.  Hence $\beta_1(M) \geq g$.  The gluing isometry $\partial M \rightarrow \partial \wwbar{M}$ (the identity) extends over $M$, thus $\mathrm{Ker}\ i_*=\mathrm{Ker}\ j_*$, and so $\mathrm{dim}\ \mathrm{Im}(i_*,-j_*)=g$.  Hence 
\[ \mathrm{H}_1(DM, \mathbb{Q}) \cong \frac{\mathrm{H}_1(M, \mathbb{Q}) \oplus \mathrm{H}_1(\wwbar{M}, \mathbb{Q})}{\mathrm{Im}((i_*,-j_*))} \]
has dimension at least $g$.

Considering the above picture gives hope that by cutting $DM$ along $\partial M$ and regluing via some isometry $\phi\co \partial M \rightarrow \partial M$ to produce a ``twisted double'' $D_{\phi}M$, one may reduce the homological contribution of $\partial M$.  For then $j=i \circ \phi$, and if $\phi_*$ moves the kernel of the inclusion off of itself, then the argument above shows that the homology of $D_{\phi}M$ will be reduced.  Below we apply this idea to a family of examples constructed by Zimmerman and Paoluzzi \cite{ZP} which build on the ``Tripos'' example of Thurston \cite{Th}.

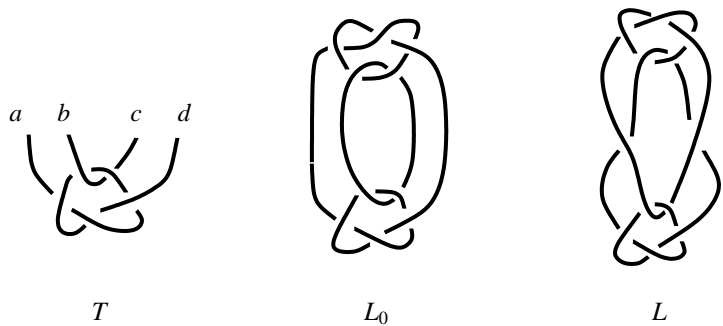
\begin{figure}[ht!]
\begin{center}
\begin{picture}(0,0)%
\includegraphics[scale=.8]{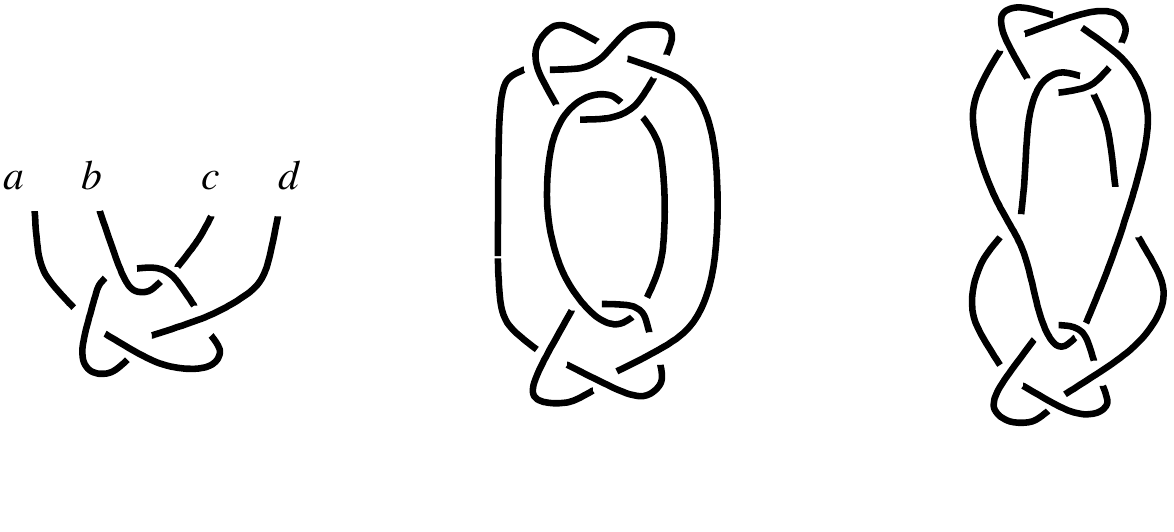}%
\end{picture}%
\setlength{\unitlength}{3315sp}%
\begingroup\makeatletter\ifx\SetFigFont\undefined%
\gdef\SetFigFont#1#2#3#4#5{%
  \reset@font\fontsize{#1}{#2pt}%
  \fontfamily{#3}\fontseries{#4}\fontshape{#5}%
  \selectfont}%
\fi\endgroup%
\begin{picture}(5354,2428)(1199,-2209)
\put(6076,-2131){\makebox(0,0)[b]{\smash{{\SetFigFont{12}{14.4}{\rmdefault}{\mddefault}{\updefault}{\color[rgb]{0,0,0}$L$}%
}}}}
\put(3961,-2131){\makebox(0,0)[b]{\smash{{\SetFigFont{12}{14.4}{\rmdefault}{\mddefault}{\updefault}{\color[rgb]{0,0,0}$L_0$}%
}}}}
\put(1891,-2131){\makebox(0,0)[b]{\smash{{\SetFigFont{12}{14.4}{\rmdefault}{\mddefault}{\updefault}{\color[rgb]{0,0,0}$T$}%
}}}}
\end{picture}%
\caption{The tangle $T$ and its double and twisted double.}
\label{links}
\end{center}
\end{figure}  

The complement in the ball of the tangle $T$ in \fullref{links} is one of the minimal volume hyperbolic manifolds with totally geodesic boundary, obtained as an identification space of a regular ideal octahedron \cite{Mi}.  We will denote it by $O_{\infty}$.  For $n \geq 3$, the orbifold $O_n$ with totally geodesic boundary consisting of the ball with cone locus $T$ of cone angle $2\pi/n$ has been explicitly described by Zimmerman and Paoluzzi \cite{ZP} as an identification space of a truncated tetrahedron.  For each $k<n$ with $(k,n)=1$, Zimmerman and Paoluzzi describe a hyperbolic manifold $M_{n,k}$ which is an $n$--fold branched cover of $O_n$.  Topologically, $M_{n,k}$ is the $n$--fold branched cover of the ball, branched over $T$, obtained as the kernel of $\langle x,y \rangle = \mathbb{Z} \oplus \mathbb{Z} \rightarrow \mathbb{Z}/n\mathbb{Z} = \langle t \rangle$ via $x \mapsto t$, $y \mapsto t^k$, where $x$ and $y$ are homology classes representing meridians of the two components of $T$.

We recall a well-known fact about isometries of spheres with four cone points:

\begin{fact} Let $S$ be a hyperbolic sphere with four cone points of equal cone angle $\alpha$, $0 \leq \alpha \leq 2\pi/3$, labeled $a$, $b$, $c$, $d$.  Each of the following permutations of the cone points may be realized by an orientation-preserving isometry:
\begin{align*}
 &(ab)(cd)& &(ac)(bd)& &(ad)(bc)&
\end{align*}
\end{fact}

Using this fact and abusing notation, let $\phi$ be the isometry $(ab)(cd)$ of $\partial O_n$, with labels as in \fullref{links}.  Doubling the tangle ball produces the link $L_0$ in \fullref{links}, and cutting along the separating $4$--punctured sphere and regluing via $\phi$ produces the link $L$, a \textit{mutant\/} of $L_0$ in the classical terminology.  Note that $L$ and all of the orbifolds $D_{\phi} O_n$ contain the mutation sphere as a totally geodesic surface, by the fact above.  $\phi$ lifts to an isometry $\tphi$ of $\partial M_{n,k}$, and the twisted double $D_{\tphi}M_{n,k}$ is the corresponding branched cover over $L$.

The homology of $D_{\tphi}M_{n,k}$ can be described using the Alexander polynomial of $L$.  The two variable Alexander polynomial of $L$ is
\[ \Delta_L(x,y) = \frac{1}{x^3}(x-1)(xy-1)(y-1)^2(x-y). \]
For the regular $\mathbb{Z}$--covering of $\mathrm{S}^3-L$ given by $x \mapsto t^k$, $y \mapsto t$, the Alexander polynomial is
\[ \Delta_L^k(t)= (t-1)\Delta(t^k,t)= \frac{1}{t^{3k-1}}(t-1)^5\nu_{k-1}(t) \nu_k(t) \nu_{k+1}(t) \]
where $\nu_k(t)=t^{k-1} + t^{k-2} + \cdots + t + 1$.  By a theorem originally due to Sumners \cite{Su} in the case of links, the first Betti number of $D_{\tphi}M_{n,k}$ is the number of roots shared by $\Delta_L^k(t)$ and $\nu_n(t)$.  Since this number is $0$ for many $n$ and $k$, we have a more precise version of \fullref{closed}.

\begin{thm1} For $n>3$ prime and $k \neq 0,1,n-1$, the manifold $D_{\tphi}M_{n,k}$ is a hyperbolic rational homology sphere containing an embedded totally geodesic surface.
\end{thm1}

The techniques used above are obviously more generally applicable.  Given any hyperbolic two-string tangle in a ball with totally geodesic boundary, one may double it to get a $2$--component hyperbolic link in $\mathbb{S}^3$ and then mutate along the separating $4$--punctured sphere by an isometry.  By the hyperbolic Dehn surgery theorem and the fact above, for large enough $n$, $(n,0)$ orbifold surgery on each component will yield a hyperbolic orbifold with a separating totally geodesic orbisurface.  Then $n$--fold manifold branched covers can be constructed as above.  One general observation about such covers follows from the following well-known fact, originally due to Conway:

\begin{fact} The one variable Alexander polynomial of a link is not altered by mutation; ie, 
\[ \Delta_{L_0}(t,t) = \Delta_L(t,t) \]
when L is obtained from $L_0$ by mutation along a $4$--punctured sphere.
\end{fact}

In our situation, this implies the following:

\begin{cor1} A $2$--component link in $\mathrm{S}^3$ which is the twisted double of a tangle has no integral homology spheres among its abelian branched covers.
\end{cor1}

\begin{proof} A link $L_0$ which is the double of a tangle has Alexander polynomial $0$.  Therefore by the fact above,
\[ \Delta_L^1(t)=(t-1)\Delta_L(t,t)=(t-1)\Delta_{L_0}(t,t)=0, \]
and so $D_{\tphi}M_{n,1}$ has positive first Betti number by Sumners' theorem.  The canonical abelian $n^2$--fold branched cover of $L$ covers $D_{\tphi}M_{n,1}$ and so also has positive first Betti number.  Since the other $n$--fold branched covers of $L$ have $n$--torsion, no branched covers of $L$ have trivial first homology.
\end{proof}

\section[Theorem 1]{\fullref{cusped}}\label{sec3}

In this section we construct hyperbolic knot complements in rational homology spheres containing closed embedded totally geodesic surfaces.  The following ``commutative diagram'' introduces the objects involved in the construction and the relationships between them.
\[ \xymatrix{ N_n \ar[rr]^{\mbox{\tiny{Dehn}}}_{\mbox{\tiny{filling}}} \ar[dd] && M_n \ar[rr]^{\mbox{\tiny{Dehn}}}_{\mbox{\tiny{filling}}} \ar[dd] && S_n \\ && && \\ N \ar[rr]^{\mbox{\tiny{orbifold}}}_{\mbox{\tiny{filling}}} && O_n && \\} \]
\fullref{cusped} may now be more precisely stated as follows.

\begin{thm1}  For each $n \geq 3$ odd, $O_n$ is a one-cusped hyperbolic orbifold containing a totally geodesic sphere with four cone points of order $n$, $M_n$ is a branched covering of $O_n$ which is a one-cusped hyperbolic manifold, and $S_n$ is a rational homology sphere.  \end{thm1}

Before beginning the proof, we give a brief sketch of the strategy.  We give an explicit polyhedral construction of a three-cusped hyperbolic manifold $N$ containing an embedded totally geodesic $4$--punctured sphere which intersects two of the cusps.  For $n \geq 3$, we give the polyhedral decomposition of the orbifold $O_n$ resulting from $n$--fold orbifold surgery on the boundary slopes of this $4$--punctured sphere.  From this it is evident that $O_n$ is hyperbolic and the sphere remains totally geodesic.  For odd $n \geq 3$, we prove that $O_n$ has a certain one-cusped $n$--fold manifold cover $M_n$ with a surgery $S_n$ which is a rational homology sphere.  This is accomplished by adapting an argument of Sakuma \cite{Sa} to relate the homology of the $n$--fold cover $N_n \rightarrow N$ corresponding to $M_n \rightarrow O_n$, to the homology of $S_n$.  $M_n$ is thus a hyperbolic knot complement in a rational homology sphere, containing the closed embedded totally geodesic surface which is a branched covering of the totally geodesic sphere with four cone points in $O_n$.

\begin{remark}  It follows from the construction that the ambient rational homology sphere $S_n$ covers an orbifold produced by $n$--fold orbifold surgery on each cusp of $N$.  Thus by the hyperbolic Dehn surgery theorem, $S_n$ is hyperbolic for $n >> 0$. \end{remark}

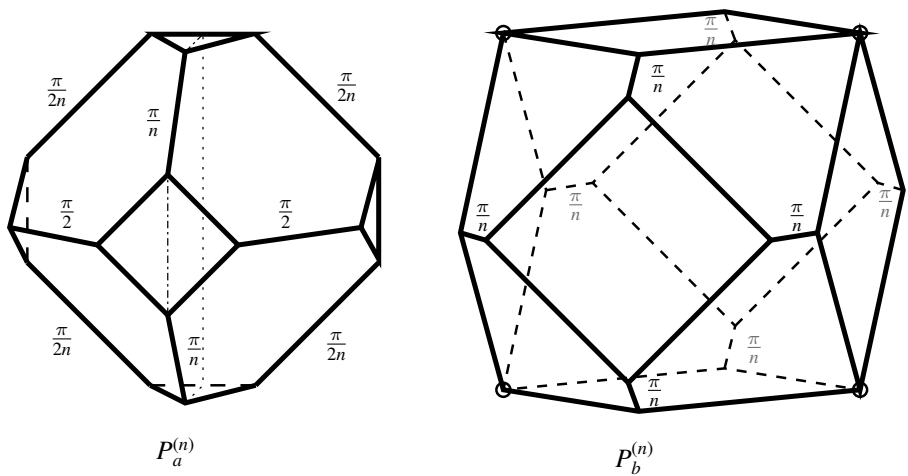
\begin{figure}[ht!]
\begin{center}
\begin{picture}(0,0)%
\includegraphics{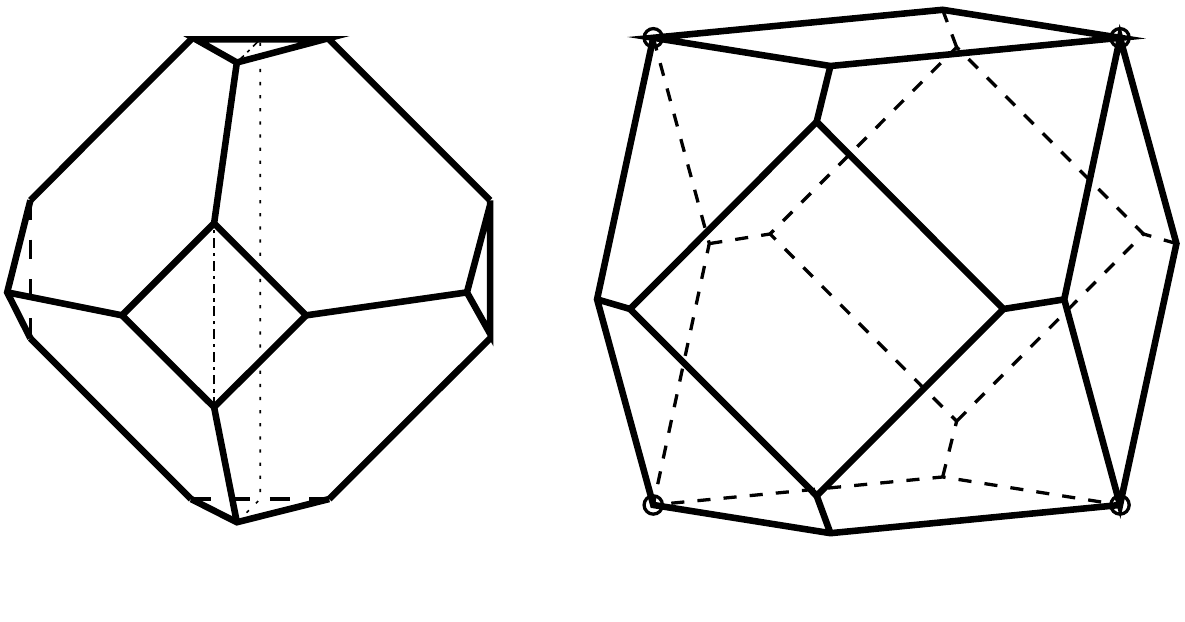}%
\end{picture}%
\setlength{\unitlength}{4144sp}%
\begingroup\makeatletter\ifx\SetFigFont\undefined%
\gdef\SetFigFont#1#2#3#4#5{%
  \reset@font\fontsize{#1}{#2pt}%
  \fontfamily{#3}\fontseries{#4}\fontshape{#5}%
  \selectfont}%
\fi\endgroup%
\begin{picture}(5412,2830)(1852,-3149)
\put(2229,-1591){\makebox(0,0)[b]{\smash{{\SetFigFont{9}{10.8}{\rmdefault}{\mddefault}{\updefault}{\color[rgb]{0,0,0}$\frac{\pi}{2}$}%
}}}}
\put(3840,-2382){\makebox(0,0)[b]{\smash{{\SetFigFont{9}{10.8}{\rmdefault}{\mddefault}{\updefault}{\color[rgb]{0,0,0}$\frac{\pi}{2n}$}%
}}}}
\put(3536,-1591){\makebox(0,0)[b]{\smash{{\SetFigFont{9}{10.8}{\rmdefault}{\mddefault}{\updefault}{\color[rgb]{0,0,0}$\frac{\pi}{2}$}%
}}}}
\put(3901,-831){\makebox(0,0)[b]{\smash{{\SetFigFont{9}{10.8}{\rmdefault}{\mddefault}{\updefault}{\color[rgb]{0,0,0}$\frac{\pi}{2n}$}%
}}}}
\put(2168,-862){\makebox(0,0)[b]{\smash{{\SetFigFont{9}{10.8}{\rmdefault}{\mddefault}{\updefault}{\color[rgb]{0,0,0}$\frac{\pi}{2n}$}%
}}}}
\put(2198,-2351){\makebox(0,0)[b]{\smash{{\SetFigFont{9}{10.8}{\rmdefault}{\mddefault}{\updefault}{\color[rgb]{0,0,0}$\frac{\pi}{2n}$}%
}}}}
\put(2989,-2412){\makebox(0,0)[b]{\smash{{\SetFigFont{9}{10.8}{\rmdefault}{\mddefault}{\updefault}{\color[rgb]{0,0,0}$\frac{\pi}{n}$}%
}}}}
\put(2746,-1044){\makebox(0,0)[b]{\smash{{\SetFigFont{9}{10.8}{\rmdefault}{\mddefault}{\updefault}{\color[rgb]{0,0,0}$\frac{\pi}{n}$}%
}}}}
\put(5693,-779){\makebox(0,0)[lb]{\smash{{\SetFigFont{7}{8.4}{\rmdefault}{\mddefault}{\updefault}{\color[rgb]{0,0,0}$\frac{\pi}{n}$}%
}}}}
\put(5201,-1548){\makebox(0,0)[lb]{\smash{{\SetFigFont{7}{8.4}{\rmdefault}{\mddefault}{\updefault}{\color[rgb]{0.447,0.431,0.427}$\frac{\pi}{n}$}%
}}}}
\put(6248,-480){\makebox(0,0)[lb]{\smash{{\SetFigFont{7}{8.4}{\rmdefault}{\mddefault}{\updefault}{\color[rgb]{0.447,0.431,0.427}\adjustlabel<-15pt,-1pt>{$\frac{\pi}{n}$}}%
}}}}
\put(7101,-1528){\makebox(0,0)[lb]{\smash{{\SetFigFont{7}{8.4}{\rmdefault}{\mddefault}{\updefault}{\color[rgb]{0.447,0.431,0.427}\adjustlabel<-3pt,0pt>{$\frac{\pi}{n}$}}%
}}}}
\put(5672,-2659){\makebox(0,0)[lb]{\smash{{\SetFigFont{7}{8.4}{\rmdefault}{\mddefault}{\updefault}{\color[rgb]{0,0,0}$\frac{\pi}{n}$}%
}}}}
\put(4646,-1613){\makebox(0,0)[lb]{\smash{{\SetFigFont{7}{8.4}{\rmdefault}{\mddefault}{\updefault}{\color[rgb]{0,0,0}$\frac{\pi}{n}$}%
}}}}
\put(6525,-1613){\makebox(0,0)[lb]{\smash{{\SetFigFont{7}{8.4}{\rmdefault}{\mddefault}{\updefault}{\color[rgb]{0,0,0}$\frac{\pi}{n}$}%
}}}}
\put(6268,-2402){\makebox(0,0)[lb]{\smash{{\SetFigFont{7}{8.4}{\rmdefault}{\mddefault}{\updefault}{\color[rgb]{0.447,0.431,0.427}$\frac{\pi}{n}$}%
}}}}
\put(2881,-3031){\makebox(0,0)[b]{\smash{{\SetFigFont{12}{14.4}{\rmdefault}{\mddefault}{\updefault}{\color[rgb]{0,0,0}$P_a^{(n)}$}%
}}}}
\put(5626,-3076){\makebox(0,0)[b]{\smash{{\SetFigFont{12}{14.4}{\rmdefault}{\mddefault}{\updefault}{\color[rgb]{0,0,0}$P_b^{(n)}$}%
}}}}
\end{picture}%
\caption{Cells for $O_n$}
\label{orbipoly}
\end{center}
\end{figure}

The proof occupies the remainder of the section.  We first discuss the orbifolds $O_n$.  For each $n$, the orbifold $O_n$ decomposes into the two polyhedra in \fullref{orbipoly}.  Realized as a hyperbolic polyhedron, $P_a^{(n)}$ is composed of two truncated tetrahedra, each of which has two opposite edges of dihedral angle $\pi/2$ and all other dihedral angles $\pi/2n$, glued along a face.  This decomposition is indicated in \fullref{orbipoly} by the lighter dashed and dotted lines.  The polyhedron $P_b^{(n)}$ has all edges with dihedral angle $\pi/2$ except for those labeled otherwise and realized as a hyperbolic polyhedron it has all combinatorial symmetries and all circled vertices at infinity.  By Andreev's theorem, polyhedra with the desired properties exist in hyperbolic space.  Certain face pairings (described below) of $P_a^{(n)}$ yield a compact hyperbolic orbifold with totally geodesic boundary a sphere with four cone points of cone angle $2\pi/n$.  Faces of $P_b^{(n)}$ may be glued to give a one-cusped hyperbolic orbifold with a torus cusp and totally geodesic boundary isometric to the boundary of the gluing of $P_a^{(n)}$.  $O_n$ is formed by gluing these orbifolds along their boundaries.

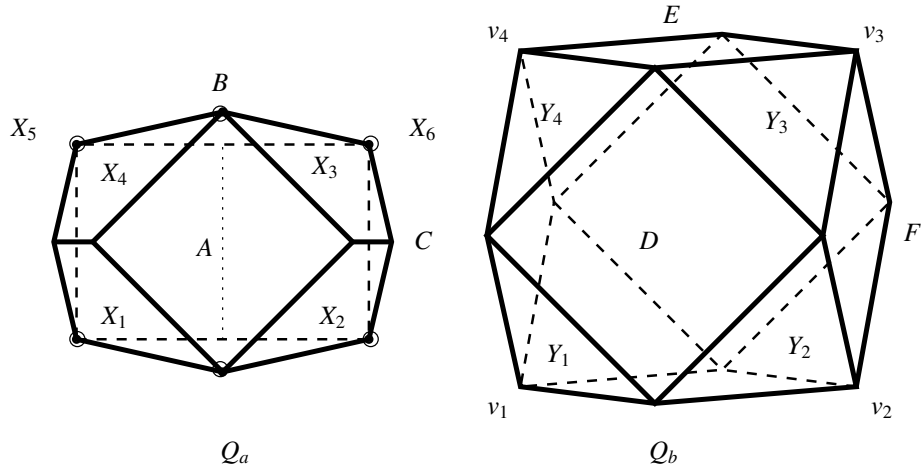
\begin{figure}[ht!]
\begin{center}
\begin{picture}(0,0)%
\includegraphics{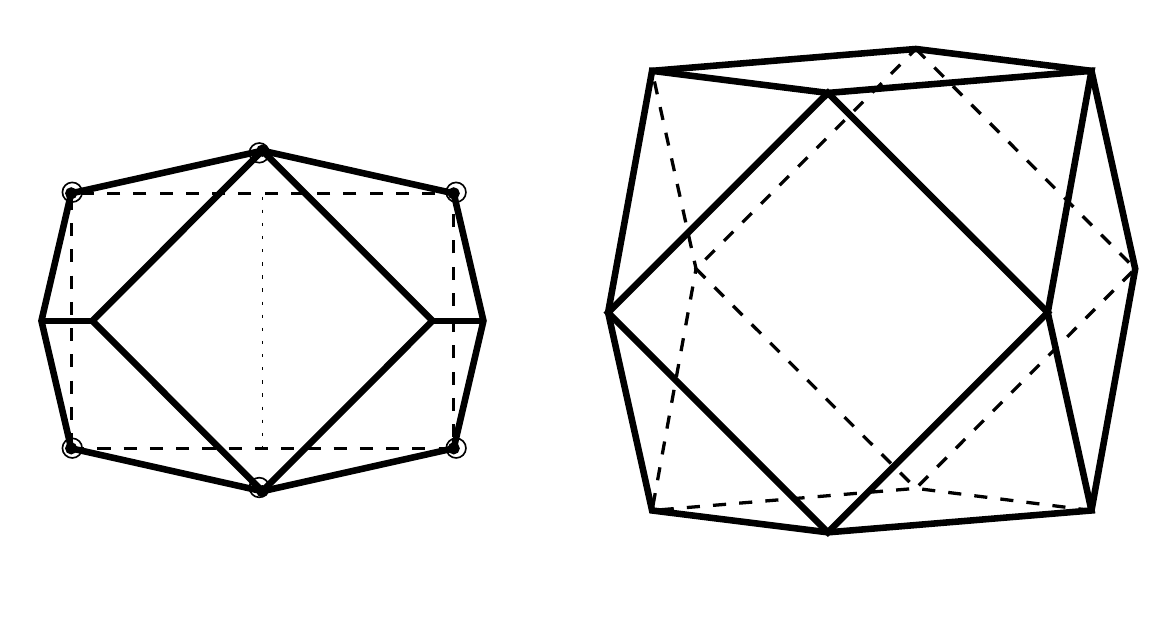}%
\end{picture}%
\setlength{\unitlength}{4144sp}%
\begingroup\makeatletter\ifx\SetFigFont\undefined%
\gdef\SetFigFont#1#2#3#4#5{%
  \reset@font\fontsize{#1}{#2pt}%
  \fontfamily{#3}\fontseries{#4}\fontshape{#5}%
  \selectfont}%
\fi\endgroup%
\begin{picture}(5340,2830)(1381,-2249)
\put(1936,-1366){\makebox(0,0)[b]{\smash{{\SetFigFont{10}{12.0}{\rmdefault}{\mddefault}{\updefault}{\color[rgb]{0,0,0}$X_1$}%
}}}}
\put(3241,-1366){\makebox(0,0)[b]{\smash{{\SetFigFont{10}{12.0}{\rmdefault}{\mddefault}{\updefault}{\color[rgb]{0,0,0}$X_2$}%
}}}}
\put(3196,-466){\makebox(0,0)[b]{\smash{{\SetFigFont{10}{12.0}{\rmdefault}{\mddefault}{\updefault}{\color[rgb]{0,0,0}$X_3$}%
}}}}
\put(1936,-511){\makebox(0,0)[b]{\smash{{\SetFigFont{10}{12.0}{\rmdefault}{\mddefault}{\updefault}{\color[rgb]{0,0,0}$X_4$}%
}}}}
\put(1396,-241){\makebox(0,0)[b]{\smash{{\SetFigFont{10}{12.0}{\rmdefault}{\mddefault}{\updefault}{\color[rgb]{0,0,0}$X_5$}%
}}}}
\put(3781,-241){\makebox(0,0)[b]{\smash{{\SetFigFont{10}{12.0}{\rmdefault}{\mddefault}{\updefault}{\color[rgb]{0,0,0}$X_6$}%
}}}}
\put(2476,-961){\makebox(0,0)[b]{\smash{{\SetFigFont{12}{14.4}{\rmdefault}{\mddefault}{\updefault}{\color[rgb]{0,0,0}$A$}%
}}}}
\put(2566, 29){\makebox(0,0)[b]{\smash{{\SetFigFont{12}{14.4}{\rmdefault}{\mddefault}{\updefault}{\color[rgb]{0,0,0}$B$}%
}}}}
\put(4591,-1591){\makebox(0,0)[b]{\smash{{\SetFigFont{10}{12.0}{\rmdefault}{\mddefault}{\updefault}{\color[rgb]{0,0,0}$Y_1$}%
}}}}
\put(6031,-1546){\makebox(0,0)[b]{\smash{{\SetFigFont{10}{12.0}{\rmdefault}{\mddefault}{\updefault}{\color[rgb]{0,0,0}$Y_2$}%
}}}}
\put(5896,-196){\makebox(0,0)[b]{\smash{{\SetFigFont{10}{12.0}{\rmdefault}{\mddefault}{\updefault}{\color[rgb]{0,0,0}$Y_3$}%
}}}}
\put(4546,-151){\makebox(0,0)[b]{\smash{{\SetFigFont{10}{12.0}{\rmdefault}{\mddefault}{\updefault}{\color[rgb]{0,0,0}$Y_4$}%
}}}}
\put(5131,-916){\makebox(0,0)[b]{\smash{{\SetFigFont{12}{14.4}{\rmdefault}{\mddefault}{\updefault}{\color[rgb]{0,0,0}$D$}%
}}}}
\put(5266,434){\makebox(0,0)[b]{\smash{{\SetFigFont{12}{14.4}{\rmdefault}{\mddefault}{\updefault}{\color[rgb]{0,0,0}$E$}%
}}}}
\put(3781,-916){\makebox(0,0)[b]{\smash{{\SetFigFont{12}{14.4}{\rmdefault}{\mddefault}{\updefault}{\color[rgb]{0,0,0}$C$}%
}}}}
\put(5221,-2176){\makebox(0,0)[b]{\smash{{\SetFigFont{12}{14.4}{\rmdefault}{\mddefault}{\updefault}{\color[rgb]{0,0,0}$Q_b$}%
}}}}
\put(2656,-2176){\makebox(0,0)[b]{\smash{{\SetFigFont{12}{14.4}{\rmdefault}{\mddefault}{\updefault}{\color[rgb]{0,0,0}$Q_a$}%
}}}}
\put(6706,-871){\makebox(0,0)[b]{\smash{{\SetFigFont{12}{14.4}{\rmdefault}{\mddefault}{\updefault}{\color[rgb]{0,0,0}$F$}%
}}}}
\put(4231,-1906){\makebox(0,0)[b]{\smash{{\SetFigFont{10}{12.0}{\rmdefault}{\mddefault}{\updefault}{\color[rgb]{0,0,0}$v_1$}%
}}}}
\put(6481,344){\makebox(0,0)[b]{\smash{{\SetFigFont{10}{12.0}{\rmdefault}{\mddefault}{\updefault}{\color[rgb]{0,0,0}$v_3$}%
}}}}
\put(4231,344){\makebox(0,0)[b]{\smash{{\SetFigFont{10}{12.0}{\rmdefault}{\mddefault}{\updefault}{\color[rgb]{0,0,0}$v_4$}%
}}}}
\put(6526,-1906){\makebox(0,0)[b]{\smash{{\SetFigFont{10}{12.0}{\rmdefault}{\mddefault}{\updefault}{\color[rgb]{0,0,0}$v_2$}%
}}}}
\end{picture}%
\caption{Cells for $N$}
\label{idealpoly}
\end{center}
\end{figure}

The geometric limit of the $O_n$ as $n \rightarrow \infty$ is $N$, a $3$--cusped manifold which decomposes into the two polyhedra in \fullref{idealpoly}.  As above, realized as a convex polyhedron in hyperbolic space $Q_a$ has all circled vertices at infinity.  The edge of $Q_a$ connecting face $A$ to face $C$ is finite length, as is the corresponding edge on the opposite vertex of $A$; all others are ideal or half-ideal and all have dihedral angle $\pi/2$.  $Q_a$ has a reflective involution of order $2$ corresponding to the involution of $P_a^{(n)}$ interchanging the two truncated tetrahedra.  The fixed set of this involution on the back face is shown as a dotted line, and notationally we regard $Q_a$ as having an edge there with dihedral angle $\pi$, splitting the back face into two faces $X_5$ and $X_6$.  $Q_b$ is the regular all-right hyperbolic ideal cuboctahedron.

Another remark on notation: the face opposite a face labeled with only a letter should be interpreted as being labeled with that letter ``prime''.  For instance, the leftmost triangular face of $Q_a$ has label $C'$.  Also, each ``back'' triangular face of $Q_b$ takes the label of the face with which it shares a vertex.  For example, the lower left back triangular face is $Y_1'$.

\begin{figure}[ht!]
\begin{center}
\begin{picture}(0,0)%
\includegraphics{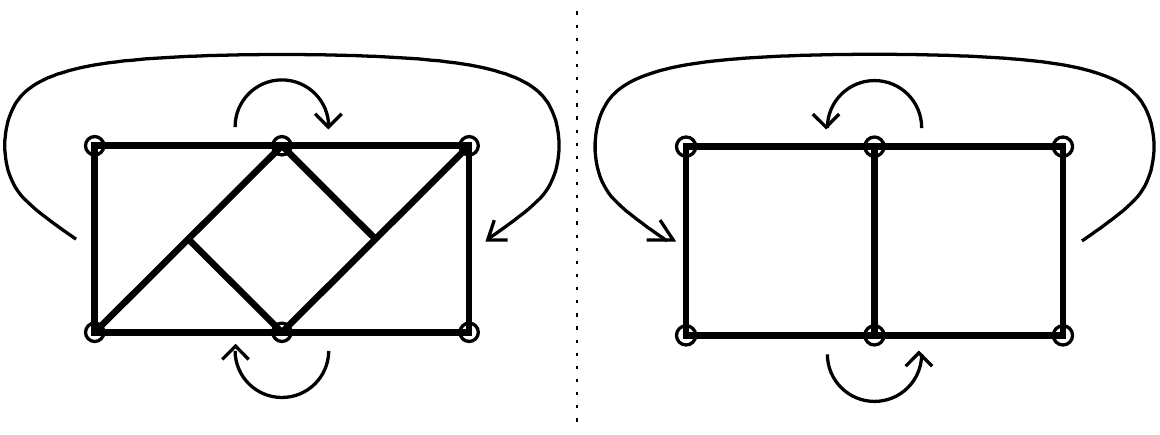}%
\end{picture}%
\setlength{\unitlength}{4144sp}%
\begingroup\makeatletter\ifx\SetFigFont\undefined%
\gdef\SetFigFont#1#2#3#4#5{%
  \reset@font\fontsize{#1}{#2pt}%
  \fontfamily{#3}\fontseries{#4}\fontshape{#5}%
  \selectfont}%
\fi\endgroup%
\begin{picture}(5298,1983)(917,-1828)
\put(2206,-1006){\makebox(0,0)[b]{\smash{{\SetFigFont{12}{14.4}{\rmdefault}{\mddefault}{\updefault}{\color[rgb]{0,0,0}$A$}%
}}}}
\put(1756,-1276){\makebox(0,0)[b]{\smash{{\SetFigFont{10}{12.0}{\rmdefault}{\mddefault}{\updefault}{\color[rgb]{0,0,0}$r^{-1}(C)$}%
}}}}
\put(2701,-1276){\makebox(0,0)[b]{\smash{{\SetFigFont{10}{12.0}{\rmdefault}{\mddefault}{\updefault}{\color[rgb]{0,0,0}$t^{-1}(B)$}%
}}}}
\put(4456,-1006){\makebox(0,0)[b]{\smash{{\SetFigFont{12}{14.4}{\rmdefault}{\mddefault}{\updefault}{\color[rgb]{0,0,0}$f_4^{-1}(D')$}%
}}}}
\put(5356,-1006){\makebox(0,0)[b]{\smash{{\SetFigFont{12}{14.4}{\rmdefault}{\mddefault}{\updefault}{\color[rgb]{0,0,0}$D$}%
}}}}
\put(5446,-376){\makebox(0,0)[b]{\smash{{\SetFigFont{12}{14.4}{\rmdefault}{\mddefault}{\updefault}{\color[rgb]{0,0,0}$f_4^{-1}f_3$}%
}}}}
\put(4231,-16){\makebox(0,0)[b]{\smash{{\SetFigFont{12}{14.4}{\rmdefault}{\mddefault}{\updefault}{\color[rgb]{0,0,0}$f_4^{-1}f_2$}%
}}}}
\put(4366,-1636){\makebox(0,0)[b]{\smash{{\SetFigFont{12}{14.4}{\rmdefault}{\mddefault}{\updefault}{\color[rgb]{0,0,0}$f_1^{-1}f_4$}%
}}}}
\put(3016,-16){\makebox(0,0)[b]{\smash{{\SetFigFont{12}{14.4}{\rmdefault}{\mddefault}{\updefault}{\color[rgb]{0,0,0}$t^{-2}s^{-2}$}%
}}}}
\put(2746,-1636){\makebox(0,0)[b]{\smash{{\SetFigFont{12}{14.4}{\rmdefault}{\mddefault}{\updefault}{\color[rgb]{0,0,0}$sts^{-1}t$}%
}}}}
\put(1576,-376){\makebox(0,0)[b]{\smash{{\SetFigFont{12}{14.4}{\rmdefault}{\mddefault}{\updefault}{\color[rgb]{0,0,0}$t^{-1}s^{-1}ts^{-1}$}%
}}}}
\put(1666,-736){\makebox(0,0)[b]{\smash{{\SetFigFont{10}{12.0}{\rmdefault}{\mddefault}{\updefault}{\color[rgb]{0,0,0}$s(B')$}%
}}}}
\put(2656,-691){\makebox(0,0)[b]{\smash{{\SetFigFont{10}{12.0}{\rmdefault}{\mddefault}{\updefault}{\color[rgb]{0,0,0}$r(C')$}%
}}}}
\end{picture}%
\caption{Totally geodesic faces of $N_a$ and $N_b$}
\label{faces}
\end{center}
\end{figure}
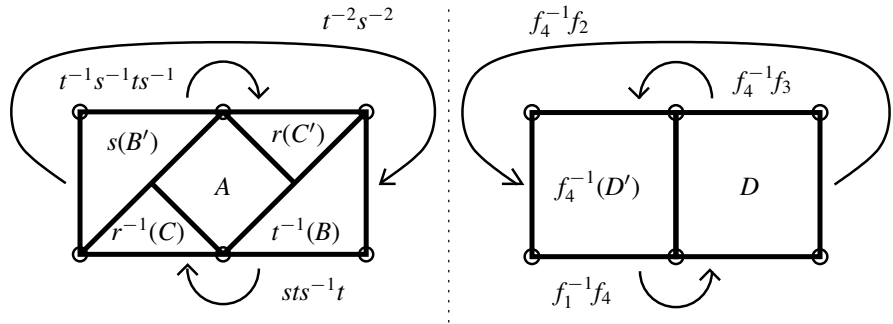

We first consider face pairings of $Q_a$ producing a manifold $N_a$ with two annulus cusps and totally geodesic boundary.  Let $r$, $s$ and $t$ be isometries realizing face pairings $X_1 \mapsto X_3$, $X_6 \mapsto X_4$ and $X_2 \mapsto X_5$, respectively.  Poincar\'{e}'s polyhedron theorem yields a presentation
\[ \langle\ r,s,t\ \,|\,\ rst=1\ \rangle \]
for the group generated by $r$, $s$ and $t$.  Note that this group is free on two generators, say $s$ and $t$, where by the relation $r=t^{-1}s^{-1}$.  Choose as the ``boundary subgroup'' (among all possible conjugates) the subgroup fixing the hyperbolic plane through the face $A$.  A fundamental polyhedron for this group and its face-pairing isometries are in \fullref{faces}.  Note that the boundary is a $4$--punctured sphere, and two of the three generators listed are the parabolics $t^{-1}s^{-1}ts^{-1}$ and $sts^{-1}t$, which generate the two annulus cusp subgroups of $\langle\,s,t\,\rangle$.

We now consider $Q_b$ and the $3$--cusped quotient manifold $N_b$.  For $i \in \{1,2,3,4\}$, let $f_i$ be the isometry pairing the face $Y_i \rightarrow Y_{i+1}'$ so that $v_i \mapsto v_{i+1}$.  Let $g_1$ be the \textit{hyperbolic\/} isometry (that is, without twisting) sending $E \rightarrow E'$ and $g_2$ the hyperbolic isometry sending $F \rightarrow F'$.  The polyhedron theorem gives presentation  
\begin{align*}
\langle f_1, f_2, f_3, f_4, g_1, g_2 \,|\, & f_1g_2f_2^{-1}g_1^{-1}=1, \\
   & f_2^{-1}g_2^{-1}f_3g_1^{-1}=1, \\ 
   & f_3g_2^{-1}f_4^{-1}g_1=1, \\ 
   & f_4^{-1}g_2f_1g_1=1 \rangle \end{align*}
for the group generated by the face pairings.  The first three generators and relations may be eliminated from this presentation using Nielsen--Schreier transformations, yielding a presentation
\[\langle\ f_4, g_1, g_2\ \,|\,\ f_4^{-1}[g_2,g_1]f_4[g_2,g_1^{-1}] = 1\ \rangle \]
(our commutator convention is $[x,y]=xyx^{-1}y^{-1}$), where the first three relations yield
\[ f_1=g_1g_2^{-1}g_1^{-1}f_4g_2g_1^{-1}g_2^{-1},\quad f_2=g_2^{-1}g_1^{-1}f_4g_2g_1^{-1}\quad \text{and}\quad f_3=g_1^{-1}f_4g_2. \]
The second presentation makes clear that the homology of $N_b$ is free of rank $3$, since each generator has exponent sum $0$ in the relation.  Faces $D$ and $D'$ make up the totally geodesic boundary of $N_b$.  \fullref{faces} shows a fundamental polyhedron for the boundary subgroup fixing $D$, together with the face pairings generating the boundary subgroup.

\begin{figure}[ht!]
\begin{center}
\begin{picture}(0,0)%
\includegraphics{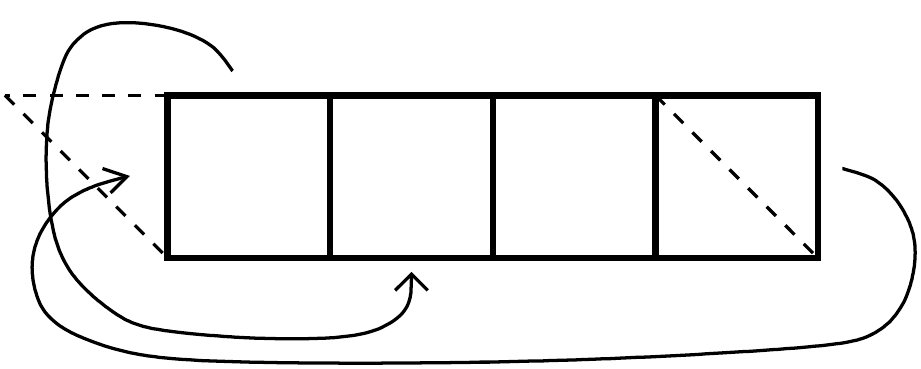}%
\end{picture}%
\setlength{\unitlength}{4144sp}%
\begingroup\makeatletter\ifx\SetFigFont\undefined%
\gdef\SetFigFont#1#2#3#4#5{%
  \reset@font\fontsize{#1}{#2pt}%
  \fontfamily{#3}\fontseries{#4}\fontshape{#5}%
  \selectfont}%
\fi\endgroup%
\begin{picture}(4206,1668)(2221,-3914)
\put(3358,-2927){\makebox(0,0)[b]{\smash{{\SetFigFont{9}{10.8}{\rmdefault}{\mddefault}{\updefault}{\color[rgb]{0,0,0}$v_1$}%
}}}}
\put(4102,-3225){\makebox(0,0)[b]{\smash{{\SetFigFont{9}{10.8}{\rmdefault}{\mddefault}{\updefault}{\color[rgb]{0,0,0}$f_1^{-1}(v_2)$}%
}}}}
\put(4845,-2927){\makebox(0,0)[b]{\smash{{\SetFigFont{9}{10.8}{\rmdefault}{\mddefault}{\updefault}{\color[rgb]{0,0,0}\adjustlabel<-25pt,22pt>{$f_1^{-1}f_2^{-1}(v_3)$}}%
}}}}
\put(5217,-3225){\makebox(0,0)[lb]{\smash{{\SetFigFont{9}{10.8}{\rmdefault}{\mddefault}{\updefault}{\color[rgb]{0,0,0}\adjustlabel<-5pt,-22pt>{$f_1^{-1}f_2^{-1}f_3^{-1}(v_4)$}}%
}}}}
\put(3284,-2370){\makebox(0,0)[b]{\smash{{\SetFigFont{10}{12.0}{\rmdefault}{\mddefault}{\updefault}{\color[rgb]{0,0,0}$bob$}%
}}}}
\put(5143,-3745){\makebox(0,0)[b]{\smash{{\SetFigFont{10}{12.0}{\rmdefault}{\mddefault}{\updefault}{\color[rgb]{0,0,0}\adjustlabel<-10pt,-5pt>{$rita$}}%
}}}}
\end{picture}%
\caption{Closed cusp of $N_b$}
\label{closedcusp}
\end{center}
\end{figure}
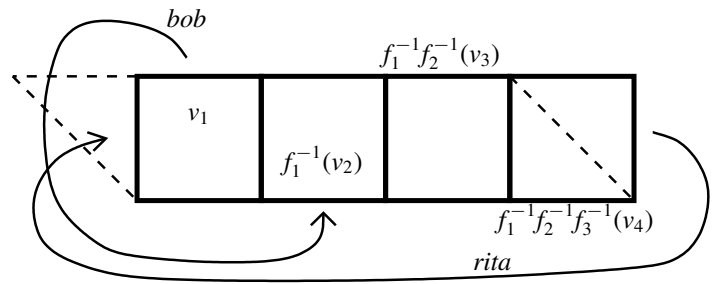

$N_b$ has two annulus cusps, each with two boundary components on the totally geodesic boundary, and one torus cusp.  A fundamental domain for the torus cusp in the horosphere centered at $v_1$ is shown in \fullref{closedcusp}, together with face pairing isometries generating the rank--$2$ parabolic subgroup fixing $v_1$.  The generators shown are
\[ bob=(f_4g_1^{-1})^2f_4g_2g_1^{-1}g_2^{-1} \quad\text{and}\quad rita=(f_4g_1^{-1})^3f_4g_2g_1^{-1}g_2^{-1}. \]
Note that $(bob)^{-4}(rita)^3$ is trivial in homology.  This and $rita \cdot (bob)^{-1}=f_4 g_1^{-1}$ together generate the cusp subgroup fixing $v_1$.  For later convenience, we now switch to the conjugate of this subgroup by $f_4^{-1}$, fixing $v_4$ and refer to the conjugated elements $m=f_4^{-1}(f_4g_1^{-1})f_4=g_1^{-1}f_4$ and $l=f_4^{-1}((bob)^{-4}(rita)^3)f_4$ as a ``meridian-longitude'' generating set for the closed cusp of $N_b$.

The totally geodesic $4$--punctured spheres on the boundaries of $N_a$ and $N_b$ are each the double of a regular ideal rectangle, and we construct $N$ by gluing $N_a$ to $N_b$ along them.  Let us therefore assume that the polyhedra in \fullref{idealpoly} are realized in hyperbolic space in such a way that face $A$ of $Q_a$ and face $D$ of $Q_b$ are in the same hyperbolic plane, with $Q_a$ and $Q_b$ in opposite half-spaces.  Further arrange so that the polyhedra are aligned in the way suggested by folding the page containing \fullref{faces} along the dotted line down the center of the figure.  With this arrangement, Maskit's combination theorem gives a presentation for the amalgamated group: \begin{align*}
 \langle\ f_4, g_1, g_2, s, t  \,|\, & 
 f_4^{-1}[g_2,g_1]f_4[g_2,g_1^{-1}] = 1, \\ 
 & t^{-2}s^{-2}=f_4^{-1}g_2^{-1}g_1^{-1}f_4g_2g_1^{-1}, \\ 
 & sts^{-1}t=g_2g_1g_2^{-1}f_4^{-1}g_1g_2g_1^{-1}f_4, \\ 
 & t^{-1}s^{-1}ts^{-1}=f_4^{-1}g_1^{-1}f_4g_2\ \rangle \end{align*}
The first relation comes from $N_b$ and the others come from setting the boundary face pairings equal to each other.  Observe that the last relation can be solved for $g_2$. Using Nielsen--Schreier transformations to eliminate $g_2$ and the last relation results in the presentation:
\begin{align*}
  \langle\ f_4, g_1, s, t\ \,|\, 
 &f_4^{-2}g_1[f_4t^{-1}s^{-1}ts^{-1},g_1]f_4[f_4t^{-1}s^{-1}ts^{-1},g_1^{-1}]g_1^{-2}f_4g_1 = 1  \\ 
 &t^{-2}s^{-2}=f_4^{-1}st^{-1}stf_4^{-1}g_1^{-1}f_4^2t^{-1}s^{-1}ts^{-1}g_1^{-1},\\ 
 &sts^{-1}t=f_4^{-1}g_1f_4t^{-1}s^{-1}ts^{-1}g_1st^{-1}stf_4^{-2}g_1f_4t^{-1}s^{-1}ts^{-1}g_1^{-1}f_4\ \rangle\end{align*}
Replace $g_1$ with the meridian generator $m=\smash{g_1^{-1}}f_4$ of the closed cusp of $N_b$ and add generators $m_1=f_4^{-1}mf_4$ and $m_2=st^{-1}stm_1t^{-1}s^{-1}ts^{-1}$, each conjugate to $m$, yielding:
\begin{align*}
 \langle\ f_4,m,m_1,m_2,s,t\ \,|\, & 
 m_1=f_4^{-1}mf_4,\ m_2=st^{-1}stm_1t^{-1}s^{-1}ts^{-1}  \\
 &m_1^{-1}t^{-1}s^{-1}ts^{-1}f_4m^{-1}m_2mf_4^{-1}st^{-1}stm_1m^{-1}=1  \\
 &s^2t^2f_4^{-1}m_2mf_4^{-1}=1  \\
 &t^{-1}st^{-1}s^{-1}m^{-1}f_4t^{-1}s^{-1}ts^{-1}f_4m^{-1}m_2^{-1}m=1\ \rangle 
\end{align*}
Note that after abelianizing, each of the last two relations expresses $f_4^2=m^2s^2t^2$, since $m_1$ and $m_2$ are conjugate to $m$ and therefore identical in homology.  In light of this, we replace $f_4$ by $u=t^{-1}s^{-1}f_4m^{-1}$, which has order $2$ in homology.  This yields the presentation:\begin{align}    
 \nonumber \langle\ m,m_1,m_2,s,t,u\ \,|\, \\
 &m_1^{-1}m^{-1}u^{-1}t^{-1}s^{-1}mstum=1 \\
 &m_2^{-1}st^{-1}stm_1t^{-1}s^{-1}ts^{-1}=1 \\
 &m_1^{-1}t^{-1}s^{-1}t^{2}um_2u^{-1}t^{-2}stm_1m^{-1}=1 \\
 &s^2t^2m^{-1}u^{-1}t^{-1}s^{-1}m_2u^{-1}t^{-1}s^{-1}=1 \\
 &t^{-1}st^{-1}s^{-1}m^{-1}stumt^{-1}s^{-1}t^2um_2^{-1}m=1\ \rangle 
\end{align}
Let $R_i$ denote the relation labeled $(i)$ in the presentation above.  In the abelianization, $R_1$ sets $m_1=m$, $R_2$ sets $m_2=m_1$, $R_3$ disappears, and the last two relations set $u^2=1$.  Therefore
\[ \mathrm{H}_1(N) \cong \mathbb{Z}^3 \oplus \mathbb{Z}/2\mathbb{Z} = \langle m \rangle \oplus \langle s \rangle \oplus \langle t \rangle \oplus \langle u \rangle. \]
(In this paper we will generally blur the distinction between elements of $\pi_1$ and their homology classes.)

The boundary slopes of the totally geodesic $4$--punctured sphere coming from $\partial N_a$ and $\partial N_b$ are represented in $\pi_1(N)$ by $t^{-1}s^{-1}ts^{-1}$ and $sts^{-1}t$.  Let $O_n$ be the finite volume hyperbolic orbifold produced by performing face identifications on $P^{(n)}_a$ and $P^{(n)}_b$ corresponding to those on $Q_a$ and $Q_b$.  $O_n$ is geometrically produced by $n$--fold orbifold filling on each of the above boundary slopes of $N$.  Appealing to the polyhedral decomposition, we see that the separating $4$--punctured sphere remains totally geodesic, becoming a sphere with four cone points of order $n$.  Our knots in rational homology spheres are certain manifold covers of the $O_n$.  In order to understand the homology of these manifold covers, we compute the homology of the corresponding abelian covers of $N$.

Let $p\co\widetilde{N} \rightarrow N$ be the maximal free abelian cover; that is, $\widetilde{N}$ is the cover corresponding to the kernel of the map $\pi_1(N) \rightarrow \mathrm{H}_1(N) \rightarrow \mathbb{Z}^3 =\langle x,y,z \rangle$ given by:
\begin{align*}
& m \mapsto x & & s \mapsto y & & t \mapsto z & & u \mapsto 1 &
\end{align*}
Let $X$ be a standard presentation $2$--complex for $\pi_1(N)$ and $\widetilde{X}$ the $2$--complex covering $X$ corresponding to $\widetilde{N} \rightarrow N$.  Then the first homology and Alexander module of $\widetilde{X}$ are naturally isomorphic to those of $\widetilde{N}$, since $N$ is homotopy equivalent to a cell complex obtained from $X$ by adding cells of dimension three and above.  The covering group $\mathbb{Z}^3$ acts freely on the chain complex of $\widetilde{X}$, so that it is a free $\mathbb{Z}[x,x^{-1},y,y^{-1},z,z^{-1}]$--module.  Below we give a presentation matrix for the Alex\-ander module of $\widetilde{X}$:
\[ \left( \begin{array}{ccccc}
 \frac{1-yz+xyz}{x^2yz} & 0 & -1 & -\frac{y^2z^2}{x} & \frac{-1+yz+z^2}{xz^2} \\
 -\frac{1}{x} & \frac{y^2}{x} & \frac{x-1}{x} & 0 & 0 \\
 0 & -\frac{1}{x} & \frac{z}{xy} & \frac{yz}{x} & -\frac{1}{x} \\
 \frac{x-1}{x^2yz} & -\frac{(x-1)(y+z)}{xz} & \frac{x-1}{xyz} & \frac{y(x-z)}{x} & \frac{1-2x+xz}{xz^2} \\
 \frac{x-1}{x^2z} & -\frac{y(x-1)(y-1)}{xz} & \frac{(x-1)(-1+y-z)}{xyz} & \frac{y(-x+xy+xyz-yz)}{x} & \frac{x+y-2xy}{xz^2} \\
 \frac{x-1}{x^2} & 0 & -\frac{z(x-1)}{xy} & -\frac{yz(x+yz)}{x} & \frac{y+xz}{xz}
\end{array} \right) \]
The rows of the matrix above correspond to lifts of the generators for $\pi_1(N)$ sharing a basepoint, ordered as $\{\tilde{m},\wwtilde{m_1},\wwtilde{m_2},\tilde{s},\tilde{t},\tilde{u}\}$ reading from the top down.  These generate $\mathrm{C}_1(\widetilde{X})$ as a $\mathbb{Z}[x,x^{-1},y,y^{-1},z,z^{-1}]$--module.  The columns are the Fox free derivatives of the relations in terms of the generators, giving a basis for the image of $\partial \mathrm{C}_2(\widetilde{X})$.  For a generator $g$ above, let $p_g$ be the determinant of the square matrix obtained by deleting the row corresponding to $\tilde{g}$.  These polynomials are:
\[ \begin{array}{lcl}
p_m & = & -(x^{-4}z^{-2})(x-1)^2(y-1)(z-1)(y+z+4yz+y^2z+yz^2) \\
p_{m_1} & = & (x^{-4}z^{-2})(x-1)^2(y-1)(z-1)(y+z+4yz+y^2z+yz^2) \\
p_{m_2} & = & -(x^{-4}z^{-2})(x-1)^2(y-1)(z-1)(y+z+4yz+y^2z+yz^2) \\
 p_s & = & (x^{-4}z^{-2})(x-1)(y-1)^2(z-1)(y+z+4yz+y^2z+yz^2) \\
 p_t & = & -(x^{-4}z^{-2})(x-1)(y-1)(z-1)^2(y+z+4yz+y^2z+yz^2) \\
 p_u & = & 0
\end{array} \]
The Alexander polynomial of $\mathrm{H}_1(\widetilde{N})$ is the greatest common factor:
\[ \Delta(x,y,z)=(x-1)(y-1)(z-1)(y+z+4yz+y^2z+yz^2) \]
up to multiplication by an invertible element of $\mathbb{Z}[x,x^{-1},y,y^{-1},z,z^{-1}]$.

Let $N_{\infty}$ be the infinite cyclic cover of $N$ factoring through $\widetilde{N}$ given by: \begin{align*}
& m \mapsto x^2 & & s \mapsto x & & t \mapsto x & & u \mapsto 1 &
\end{align*}
Then the chain complex of $N_{\infty}$ is a $\Lambda$--module, where $\Lambda=\mathbb{Z}[x,x^{-1}]$ and specializing the above picture yields an Alexander polynomial 
\[ \begin{array}{lcl}
  \Delta_{\infty}(x) & = & (x^2-1)(x-1)^2(2x+4x^2+2x^3) \\
   & = & 2x(x-1)^3(x+1)^3. \end{array} \]
Let $N_n$ be the $n$--fold cyclic cover of $N$ factoring through $N_{\infty}$.  For $n$ odd, $N_n$ has three cusps, since $m$, $sts^{-1}t$, and $t^{-1}s^{-1}ts^{-1}$ map to $x^{\pm2}$, which generates $\mathbb{Z}/n\mathbb{Z}$.  Let $S_n$ be the closed manifold obtained by filling $N_n$ along the slopes covering $m$, $sts^{-1}t$, and $t^{-1}s^{-1}ts^{-1}$.  Theorem \ref{cusped} follows quickly from the following lemma.

\begin{lem} For odd $n \geq 3$, $S_n$ is a rational homology sphere.  \end{lem}

\begin{proof} The proof is adapted from an analogous proof of Sakuma concerning link complements in $S^3$.

The chain complex of $N_n$ is isomorphic to $C_*(N_{\infty}) \otimes (\Lambda/(x^n-1))$.  Note that $x^n - 1 = (x-1)\nu_n$, where $\nu_n(x) = x^{n-1}+x^{n-2} + \hdots + x + 1$.  Sakuma observes that the short exact sequence of coefficient modules
\[ 0 \rightarrow \mathbb{Z} \stackrel{\nu_n}{\longrightarrow} \Lambda/(x^n-1) \rightarrow \Lambda/(\nu_n) \rightarrow 0 \]
where the map on the left is multiplication by $\nu_n$, gives rise to a short exact sequence in homology
\[ 0 \rightarrow \mathrm{H}_1(N) \stackrel{tr}{\longrightarrow} \mathrm{H}_1(N_n) \rightarrow \mathrm{H}_1(N_{\infty})/\nu_n \mathrm{H}_1(N_{\infty}) \rightarrow 0 \]
where $tr$ is the transfer map, $tr(h)=h+x.h+\cdots+x^{n-1}.h$ for a homology class $h$. Define $\mathrm{H}_n=\mathrm{H}_1(N_{\infty})/\nu_n \mathrm{H}_1(N_{\infty})$.  Since the Alexander polynomial of $N_{\infty}$ does not share roots with $\nu_n$, $\mathrm{H}_n$ is a torsion $\mathbb{Z}$--module.

The lemma follows from a comparison between $\mathrm{H}_1(S_n)$ and $\mathrm{H}_n$.  The Mayer--Vietoris sequence implies that $\mathrm{H}_1(S_n)$ is obtained as the quotient of $\mathrm{H}_1(N_n)$ by the subgroup generated by transfers of the meridians.  If $N$ were a link complement in $\mathrm{S}^3$, it would immediately follow that $\mathrm{H}_n = \mathrm{H}_1(S_n)$, since the homology of a link complement is generated by meridians.  In our case we have
\[ \mathrm{H}_n = \mathrm{H}_1(N_n) / \langle\ tr(m),tr(s),tr(t),tr(u)\ \rangle, \]
whereas 
\[ \mathrm{H}_1(S_n) = \mathrm{H}_1(N_n) / \langle\ tr(m),tr(2s),tr(2t)\ \rangle. \]
However one observes that $\mathrm{H}_1(S_n) \rightarrow \mathrm{H}_n$ is an extension of degree at most 8 (since $u$ has order 2 in $\mathrm{H}_1(N)$), and so $\mathrm{H}_1(S_n)$ is also a torsion group.
\end{proof}

Let $M_n$ be the manifold obtained by filling two of the three cusps of $N_n$ along the slopes covering $sts^{-1}t$ and $t^{-1}s^{-1}ts^{-1}$.  We have geometrically described $M_n$ as a branched cover of $O_n$, produced by $n$--fold orbifold filling along $sts^{-1}t$ and $t^{-1}s^{-1}ts^{-1}$.  There is a closed totally geodesic surface in $M_n$ covering the totally geodesic sphere with four cone points in $O_n$.  A closed manifold $S_n$ is produced by filling the remaining cusp of $M_n$ along the meridian covering $m$.  Since $S_n$ is a rational homology sphere, $M_n$ is a knot complement in a rational homology sphere, and we have proven \fullref{cusped}.

\section{Further directions}

Performing ordinary Dehn filling along the three meridians of $N$ specified in the previous section yields a manifold $S$, which is easily seen to be the connected sum of two spherical manifolds.  The half arising from the truncated tetrahedra is the quotient of $\mathrm{S}^3$, regarded as the set of unit quaternions, by the subgroup $\langle i,j,k \rangle$.  The half arising from the cuboctahedron is the lens space $\mathrm{L}(4,1)$.  The manifolds $S_n$ may be regarded as $n$--fold branched covers over the three-component link $L$ in $S$ consisting of the cores of the filling tori.

Since the meridians $t^{-1}s^{-1}ts^{-1}$ and $sts^{-1}t$ represent squares of primitive elements in the homology of $N$, any cover of $S$ branched over $L$ will have nontrivial homology of order $2$ coming from the transfers of $s$ and $t$.  However, it is possible that techniques similar to those above may be used to create knot complements in integral homology spheres.  If the manifold $N$ above---in addition to its geometric properties---had trivial nonperipheral \textit{integral\/} homology, then $S$ would be an integral homology sphere.  Porti \cite{Po} has supplied a formula in terms of the Alexander polynomial for the order of the homology of a cover of an integral homology sphere branched over a link, generalizing work of Mayberry--Murasugi in the case of $S^3$ \cite{MM}.  Using this formula, the order of the homology of branched covers of $S$ could be easily checked.

In fact, the Menasco--Reid conjecture itself may be approached using a variation of these techniques.  A genus $n-1$ handlebody may be obtained as the $n$--fold branched cover of a ball over the trivial $2$--string tangle, so knot complements in the genus $n-1$ handlebody may be obtained as $n$--fold branched covers over the trivial tangle of a knot complement in the ball.  In analogy with \fullref{sec3}, allowing the complement of $T$ to play the role of $N_a$ we ask the following:

\begin{ques} Does there exist a hyperbolic $3$--manifold with one rank $2$ and two rank $1$ cusps, which is the complement of a tangle in the ball, with totally geodesic boundary isometric to the totally geodesic boundary of the complement of the tangle $T$? \end{ques}

Such a manifold would furnish an analog of the manifold $N_b$ in \fullref{sec3}.  If the glued manifold $N$ was a $2$--component link complement in $\mathrm{S}^3$, with an unknotted component intersecting the totally geodesic Conway sphere, and this sphere remained totally geodesic under the right orbifold surgery along its boundary slopes, branched covers would give a counterexample to the Menasco--Reid conjecture.  In any case, Thurston's hyperbolic Dehn surgery theorem implies that as $n \rightarrow \infty$, the resulting surfaces would have principal curvature approaching 0, furnishing new examples of the phenomenon discovered by Leininger in \cite{Le} (although unlike Leininger's examples this would not give bounded genus).

\bibliographystyle{gtart}
\bibliography{link}

\begin{thebibliography}{}
\providecommand\bibmarginpar{\leavevmode\marginpar}
\def\urlstyle#1{{\tt #1}}

\bibitem{Ad}
\textbf{C\,C Adams}, \href{http://dx.doi.org/10.1016/0040-9383(94)90017-5}
  {\emph{Toroidally alternating knots and links}}, Topology 33 (1994) 353--369
  \xox{MR}{1273788}

\bibitem{Aetal}
\textbf{C\,C Adams}, \textbf{J\,F Brock}, \textbf{J Bugbee}, \textbf{et~al},
  \href{http://dx.doi.org/10.1016/0166-8641(92)90130--R} {\emph{Almost
  alternating links}}, Topology Appl. 46 (1992) 151--165 \xox{MR}{1184114}

\bibitem{AR}
\textbf{C\,C Adams}, \textbf{A\,W Reid}, \emph{Quasi--{F}uchsian surfaces in
  hyperbolic knot complements}, J. Austral. Math. Soc. Ser. A 55 (1993)
  116--131 \xox{MR}{1231698}

\bibitem{Ha}
\textbf{A Hatcher}, \emph{Basic Topology of Three Manifolds}, online notes
\ Available at \setbox0\hbox{\makeatletter\@url
{http://www.math.cornell.edu/~hatcher/3M/3Mdownloads.html}}
\href{http://www.math.cornell.edu/~hatcher/3M/3Mdownloads.html}
{\unhbox0}

\bibitem{IO}
\textbf{K Ichihara}, \textbf{M Ozawa},
  \href{http://dx.doi.org/10.1142/S0218216500000414} {\emph{Accidental surfaces
  in knot complements}}, J. Knot Theory Ramifications 9 (2000) 725--733
  \xox{MR}{1775383}

\bibitem{Le}
\textbf{C\,J Leininger}, \href{http://dx.doi.org/10.1142/S0218216506004531}
  {\emph{Small curvature surfaces in hyperbolic 3--manifolds}}, J. Knot Theory
  Ramifications 15 (2006) 379--411 \xox{MR}{2217503}

\bibitem{Long}
\textbf{D\,D Long}, \emph{Immersions and embeddings of totally geodesic
  surfaces}, Bull. London Math. Soc. 19 (1987) 481--484 \xox{MR}{898729}

\bibitem{LP}
\textbf{M\,T Lozano}, \textbf{J\,H Przytycki}, \emph{Incompressible surfaces in
  the exterior of a closed 3--braid I: {S}urfaces with horizontal boundary
  components}, Math. Proc. Cambridge Philos. Soc. 98 (1985) 275--299
  \xox{MR}{795894}

\bibitem{Ma}
\textbf{H Matsuda}, \href{http://dx.doi.org/10.1016/S0166-8641(01)00034-7}
  {\emph{Complements of hyperbolic knots of braid index four contain no closed
  embedded totally geodesic surfaces}}, Topology Appl. 119 (2002) 1--15
  \xox{MR}{1881706}

\bibitem{MM}
\textbf{J\,P Mayberry}, \textbf{K Murasugi}, \emph{Torsion-groups of abelian
  coverings of links}, Trans. Amer. Math. Soc. 271 (1982) 143--173
  \xox{MR}{648083}

\bibitem{MeR}
\textbf{W Menasco}, \textbf{A\,W Reid}, \emph{Totally geodesic surfaces in
  hyperbolic link complements}, from: ``Topology '90 (Columbus, OH, 1990)'',
  Ohio State Univ. Math. Res. Inst. Publ. 1, de Gruyter, Berlin (1992)
  215--226 \xox{MR}{1184413}

\bibitem{Mi}
\textbf{Y Miyamoto}, \href{http://dx.doi.org/10.1016/0040-9383(94)90001-9}
  {\emph{Volumes of hyperbolic manifolds with geodesic boundary}}, Topology 33
  (1994) 613--629 \xox{MR}{1293303}

\bibitem{Oer}
\textbf{U Oertel},
  \href{http://projecteuclid.org/getRecord?id=euclid.pjm/1102710789}
  {\emph{Closed incompressible surfaces in complements of star links}}, Pacific
  J. Math. 111 (1984) 209--230 \xox{MR}{732067}

\bibitem{ZP}
\textbf{L Paoluzzi}, \textbf{B Zimmermann}, \emph{On a class of hyperbolic
  3--manifolds and groups with one defining relation}, Geom. Dedicata 60 (1996)
  113--123 \xox{MR}{1384422}

\bibitem{Po}
\textbf{J Porti}, \emph{Mayberry--Murasugi's formula for links in homology
  3--spheres} \xox{arXiv}{math/GT0306181}

\bibitem{Sa}
\textbf{M Sakuma}, \emph{The homology groups of abelian coverings of links},
  Math. Sem. Notes Kobe Univ. 7 (1979) 515--530 \xox{MR}{567241}

\bibitem{Su}
\textbf{D\,W Sumners},
  \href{http://links.jstor.org/sici?sici=0002-9939(197410)46:1%3C143:OTHOFC%3E%
2.0.CO%3B2-5} {\emph{On the homology of finite cyclic coverings of
  higher-dimensional links}}, Proc. Amer. Math. Soc. 46 (1974) 143--149
  \xox{MR}{0350747}

\bibitem{Th}
\textbf{W\,P Thurston}, \emph{The Geometry and Topology of 3--manifolds}
  (1978)\ \ Mimeographed lecture notes

\end{thebibliography}

\end{document}